\documentclass[preprint,12pt]{elsarticle}

\usepackage{amssymb}
 \usepackage{amsthm}
\usepackage{ amsfonts}
\usepackage{amsmath}

 \biboptions{sort&compress}

\journal{Journal of Approximation Theory}

\begin{document}

\begin{frontmatter}



\title{Exact values of Kolmogorov widths of classes of Poisson integrals}


\author{A.S. Serdyuk}
\ead{serdyuk@imath.kiev.ua}
\author{V.V. Bodenchuk\corref{cor1}}
\ead{bodenchuk@ukr.net}
\cortext[cor1]{Corresponding author.}

\address{Institute of Mathematics of Ukrainian National Academy of Sciences, 3, Tereshchenkivska st., 01601, Kyiv-4, Ukraine}

\begin{abstract}
We prove that the Poisson kernel $P_{q,\beta}(t)=\sum\limits_{k=1}^{\infty}q^k\cos\left(kt-\dfrac{\beta\pi}{2}\right)$, ${q\in(0,1)}$, $\beta\in\mathbb{R}$,
satisfies Kushpel's condition $C_{y,2n}$ beginning with a number $n_q$ where $n_q$ is the smallest number $n\geq9$, for which the following inequality is satisfied:
\begin{equation*}
\dfrac{43}{10(1-q)}q^{\sqrt{n}}+\dfrac{160}{57(n-\sqrt{n})}\; \dfrac{q}{(1-q)^2}\leq
\left(\dfrac{1}{2}+\dfrac{2q}{(1+q^2)(1-q)}\right)\left(\dfrac{1-q}{1+q}\right)^{\frac {4}{1-q^2}}.
\end{equation*}
As a consequence, for all $n\geq n_q$ we obtain lower bounds for Kolmogorov widths in the space $C$ of classes $C_{\beta,\infty}^q$ of Poisson integrals of functions that belong to the unit ball in the space $L_\infty$.
The obtained estimates coincide with the best uniform approximations by trigonometric polynomials for these classes.
As a result, we obtain exact values for widths of classes $C_{\beta,\infty}^q$ and show that subspaces of trigonometric polynomials of order $n-1$ are optimal for widths of dimension $2n$.
\end{abstract}


\end{frontmatter}


\section{Introduction}
\label{}




Let $L=L_1$ denote the space of $2\pi$-periodic summable functions $f$ with the norm
$\|f\|_1=\int\limits_{-\pi}^{\pi}|f(t)|dt$,
$L_\infty$ be the space of $2\pi$-periodic measurable and essentially bounded functions with the norm
$\|f\|_\infty=\mathop{\rm ess\,sup}\limits_{t\in\mathbb{R}\ }|f(t)|,$
and $C$ be the space of $2\pi$-periodic continuous functions $f$ with the norm defined by the equality $\|f\|_C=\max\limits_{t\in\mathbb{R}}|f(t)|$.

Let $\Psi_\beta(t)$ be a fixed summable kernel of the form
\begin{equation}\label{Psi_beta}
\Psi_\beta(t)=\sum\limits_{k=1}^{\infty}\psi(k)\cos\left(kt-\dfrac{\beta\pi}{2}\right),\,\psi(k)>0,\, \sum\limits_{k=1}^{\infty}\psi(k)<\infty,\, \beta\in\mathbb{R}.
\end{equation}
The notation $C_{\beta,p}^{\psi}$, $p=1, \infty$, is used for the class of $2\pi$-periodic functions $f$ that are convolutions with the kernel $\Psi_\beta$:
\begin{equation}\label{f}
f(x)=A+\dfrac{1}{\pi}\int\limits_{-\pi}^{\pi}\Psi_\beta(x-t)\varphi(t)dt=A+\left(\Psi_\beta\ast\varphi\right)(x),\; A\in\mathbb{R},
\end{equation}
where
\begin{equation*}
\|\varphi\|_p\leq1,\, \varphi\perp1.
\end{equation*}
The function $\varphi$ in \eqref{f} is called $(\psi,\beta)$-derivative of the function $f$ and denoted by $f_\beta^\psi$.
The notion of $(\psi,\beta)$-derivative was suggested by O.I.~Stepanets (see e.g. \cite[chapter 3, sections~7--8]{Stepanets_2005}).

An important particular case for kernels $\Psi_\beta(t)$ of the form \eqref{Psi_beta} with $\psi(k)=q^k$, $q\in(0,1)$, is represented by Poisson kernels $P_{q,\beta}(t)$ with parameters $q$ and $\beta$, i.e., by functions of the form
\begin{equation*}
P_{q,\beta}(t)=\sum\limits_{k=1}^{\infty}q^k\cos\left(kt-\dfrac{\beta\pi}{2}\right),\, q\in(0,1),\, \beta\in\mathbb{R}.
\end{equation*}
Functions $f$, which are representable in the form of convolution \eqref{f} with the kernel ${\Psi_\beta(t)=P_{q,\beta}(t)}$, are called Poisson integrals.
In this case, classes $C_{\beta,p}^{\psi}$ are denoted by $C_{\beta,p}^{q}$, and $(\psi,\beta)$-derivatives $f_\beta^\psi$ of a function $f\in C_{\beta,p}^{\psi}$ with $\psi(k)=q^k$ are denoted by $f_\beta^q$.

The notation $E_n(C_{\beta,p}^{\psi})_X$, where $p=1,\infty$ and $X=L$ or $C$ respectively, is used for the best approximation of the class $C_{\beta,p}^{\psi}$ in the space $X$ by the subspace $\mathcal{T}_{2n-1}$ of trigonometric polynomials $t_{n-1}$ of order $n-1$, i.e.,
\begin{equation}\label{E_n}
E_n(C_{\beta,p}^{\psi})_X=\sup_{f\in C_{\beta,p}^{\psi}}\inf_{t_{n-1}\in \mathcal{T}_{2n-1}}\|f-t_{n-1}\|_X,
\end{equation}
and $d_m(C_{\beta,p}^{\psi},X)$ denotes the Kolmogorov width of order $m$ for the class $C_{\beta,p}^{\psi}$ in the space $X$, i.e.,
\begin{equation}\label{d_m}
d_m(C_{\beta,p}^{\psi},X)=\inf_{F_m\subset X}\sup_{f\in C_{\beta,p}^{\psi}}\inf_{y\in F_m}\|f-y\|_X,
\end{equation}
where the outer $\inf$ operation is considered over all $m$-dimensional linear subspaces $F_m$ from $X$.

The problem of finding Kolmogorov widths of functional compacts in various functional spaces has a rich history, which you can read about in books \cite{Pinkus_1985,Tihomirov_1976,Kornejchuk_1987}.
In our work, we consider the problem of finding exact values of widths $d_{2n}(C_{\beta,\infty}^{q},C)$, $d_{2n-1}(C_{\beta,\infty}^{q},C)$, and $d_{2n-1}(C_{\beta,1}^q,L)$ for all natural $n$ exceeding a number that depends on $q$ only.

Note that the problem of finding exact values of widths of classes of Poisson integrals $C_{\beta,p}^{q}$, $p=1,\infty$, is considered in \cite{Kuspel_1988,Kuspel_1989,Shevaldin_1992,Nguen_1994}.
As for the widths of sets of Poisson integrals for functions from classes $H_\omega$ that are generated by the module of continuity, sharp or asymptotically sharp estimates are obtained only in certain cases (see e.g. \cite{Shevaldin_1994,Serdyuk_2011}).

For quantities \eqref{E_n} and \eqref{d_m}, we have the relation
\begin{equation}\label{d_m_E_n}
d_{2n-1}(C_{\beta,p}^{\psi},X)\leq E_n(C_{\beta,p}^{\psi})_X.
\end{equation}

As it follows from \cite{Shevaldin_1992} (see also \cite{Bushanskiy_1978}), for any  $q\in(0,1)$, $\beta\in\mathbb{R}$, and $n\in\mathbb{N}$, we have
\begin{equation*}
E_n(C_{\beta,\infty}^{q})_C=E_n(C_{\beta,1}^{q})_{L}=\|P_{q,\beta}\ast\varphi_n\|_C=
\end{equation*}
\begin{equation}\label{E_n_rivnosti}
=\dfrac{4}{\pi}\left|\sum\limits_{\nu=0}^{\infty}\dfrac{q^{(2\nu+1)n}}{2\nu+1}\sin\left((2\nu+1)\theta_n\pi-\dfrac{\beta\pi}{2}\right)\right|,
\end{equation}
where
\begin{equation}\label{varp}
\varphi_n(t)=\textnormal{sgn}\sin nt
\end{equation}
and $\theta_n=\theta_n(q,\beta)$ is the unique root of the equation
\begin{equation}\label{theta}
\sum\limits_{\nu=0}^{\infty}q^{(2\nu+1)n}\cos\left((2\nu+1)\theta_n\pi-\dfrac{\beta\pi}{2}\right)=0
\end{equation}
on $[0,1)$
(for $\beta\in\mathbb{Z}$ equalities \eqref{E_n_rivnosti} are obtained in \cite{Krein_1938,Nikolskiy_1946}).
Therefore, in order to solve the problem on exact values for the widths under consideration, it remains to prove the following lower estimates:
\begin{equation}\label{dno1}
d_{2n}(C_{\beta,\infty}^q, C)\geq\|P_{q,\beta}\ast\varphi_n\|_C,
\end{equation}
\begin{equation}\label{dno2}
d_{2n-1}(C_{\beta,1}^q, L)\geq\|P_{q,\beta}\ast\varphi_n\|_C.
\end{equation}

The problem of obtaining estimates \eqref{dno1} and \eqref{dno2} faces fundamental difficulties related to the fact that Poisson kernels $P_{q,\beta}(t)$ may increase oscillations. In particular, as it was shown in \cite[p.~1318--1319]{Kuspel_1988}, $P_{q,\beta}\not\in \text{CVD}$ for $q=1/7$ and $\beta=0$, where $\text{CVD}=\{K\in L: \nu(K\ast f)\leq\nu(f), f\in C\}$, $\nu(g)$ is the number of sign changes of function $g\in C$ on $[0,2\pi)$.
Hence, for classes of convolutions with kernels $P_{q,\beta}(t)$, we can not obtain exact lower bounds for widths by using methods and approaches suggested by A.~Pinkus \cite{Pinkus_1985}.

Up to date, estimates \eqref{dno1} and \eqref{dno2} were known in the following cases:
\begin{itemize}
\item
For any $n\in\mathbb{N}$ and $\beta\in\mathbb{R}$ with $0<q\leq q(\beta)$ where $q(\beta)=0{.}2$ for $\beta\in\mathbb{Z}$ or $q(\beta)=0{.}196881$ for $\beta\in\mathbb{R\setminus Z}$ (see \cite{Shevaldin_1992});
\item
For $\beta=2kl$, $l\in\mathbb{Z}$, any $0<q<1$, and all numbers $n$ exceeding a number $n_*$ that depends on $q$ (in this case, it was proved that $n_*$ does exist but no way was suggested to find it constructively) \cite{Nguen_1994}.
\end{itemize}

All results in these cases were obtained on the base of Kushpel's method \cite{Kuspel_1988} of finding lower bounds for widths of classes of convolutions with generator kernels $\Psi_\beta$ under so called condition $C_{y,2n}$.
In the present paper, we also use this approach.

Let us recall some definitions and formulate known results that will be used in order to represent the results of our work.

Let $\Delta_{2n}=\{0=x_0<x_1<\dots<x_{2n}=2\pi\}$, $x_k=k\pi/n$, denote a partition of the interval $[0,2\pi]$. Consider the function
\begin{equation}\label{Ps}
\Psi_{\beta,1}(t)=(\Psi_{\beta}\ast B_{1})(t)=\sum\limits_{k=1}^{\infty}\frac{\psi(k)}{k}\cos\left(kt-\dfrac{(\beta+1)\pi}{2}\right),
\end{equation}
where $B_{1}=\sum\limits_{k=1}^{\infty}k^{-1}\sin kt$ is the Bernoulli kernel.
Let $S\Psi_{\beta,1}(\Delta_{2n})$ denote the space $S\Psi_{\beta,1}(\cdot)$ of $SK$-splines over the partition $\Delta_{2n}$, i.e., the set of functions
\begin{equation}\label{SK}
S\Psi_{\beta,1}(\cdot)=\alpha_{0}+\sum\limits_{k=1}^{2n} \alpha_k \Psi_{\beta,1}(\cdot-x_k),\;\sum\limits_{k=1}^{2n} \alpha_k=0,
\end{equation}
\begin{equation*}
\alpha_k\in\mathbb{R},\; k= 0, 1, \dots, 2n.
\end{equation*}
A fundamental $SK$-spline is a function $\overline{S\Psi}_{\beta,1}(\cdot)=\overline{S\Psi}_{\beta,1}(y, \cdot)$ of the form \eqref{SK} that satisfies the relation
\begin{equation*}
\overline{S\Psi}_{\beta,1}(y, y_k)=\delta_{0,k}=
\begin{cases}
0, & k=\overline{1,2n-1},    \\
1, & k=0,
\end{cases}
\end{equation*}
where $y_k=x_k+y$, $x_k=k\pi/n$, $y\in[0,\dfrac{\pi}{n})$.
The spline $\overline{S\Psi}_{\beta,1}(y,\cdot)$ generates a system of fundamental splines of the form $\overline{S\Psi}_{\beta,1}(y,\cdot-x_k)$, $k=\overline{0, 2n-1}$, and this system represents a basis in the space $S\Psi_{\beta,1}(\Delta_{2n})$.
Necessary and sufficient conditions of existence and uniqueness for the fundamental spline $\overline{S\Psi}_{\beta,1}(y,\cdot)$, which depends on the relation between $y$ (it is a displacement of interpolation nodes) and parameters $\psi$ and $\beta$ of the generator kernel $\Psi_{\beta,1}$, were studied in \cite{Kuspel_1988,Shevaldin_1989,Shevaldin_1993,Stepanets_1994,Serdyuk_1998,Serdyuk_1999_4}.

Since, due to the definition of $(\psi,\beta)$-derivative, for the kernel $\Psi_{\beta,1}$ we have the equality \begin{equation}\label{e15}
(\Psi_{\beta,1}(\cdot))_\beta^\psi=B_{1}(\cdot),
\end{equation}
as a consequence of \eqref{SK} we also have the equality
\begin{equation}\label{e16}
(S\Psi_{\beta,1}(\cdot))_\beta^\psi=\sum\limits_{k=1}^{2n}\alpha_kB_{1}(\cdot-x_k),\;\sum\limits_{k=1}^{2n} \alpha_k=0.
\end{equation}
Equalities in \eqref{e15} and \eqref{e16} are understood as the equality of two functions in $L$ (i.e., almost everywhere).
Due to Lemma 2.3.4 from \cite[p.~76]{Kornejchuk_1987}, the function in the right-hand side of equality \eqref{e16} is constant on each interval $(x_k, x_{k+1})$.
Therefore, among $(\psi, \beta)$-derivatives of any spline of the form \eqref{SK}, and hence, for the fundamental spline $\overline{S\Psi}_{\beta,1}(\cdot)$ as well, one can find a function which is constant on each interval $(x_k, x_{k+1})$.
In what follows, $(\overline{S\Psi}_{\beta,1}(\cdot))_\beta^\psi$ will denote such a function only.

\textbf{Definition (A.K.~Kushpel).}
\textit{For a real number $y$ and partition $\Delta_{2n}$, we say that a kernel $\Psi_\beta(\cdot)$ of the form \eqref{Psi_beta} satisfies the condition $C_{y,2n}$ (and we write $\Psi_\beta\in C_{y,2n}$) if there exists just one fundamental spline $\overline{S\Psi}_{\beta,1}(y,\cdot)$ for this kernel and the following equalities are satisfied}:
\begin{equation*}
\textrm{sgn}(\overline{S\Psi}_{\beta,1}(y,t_k))_\beta^\psi=(-1)^k\varepsilon e_k,\, k=\overline{0,2n-1},
\end{equation*}
\textit{where $t_k=(x_{k}+x_{k+1})/2,$ $e_k$ is equal to either 0 or 1, and $\varepsilon$ takes values $\pm1$ and does not depend on $k$.}

The following theorem provides a possibility to find lower estimates for Kolmogorov widths of convolution classes which are generated by kernels under the condition $C_{y,2n}$.

\textbf{Theorem 1 (A.K.~Kushpel \cite{Kuspel_1988,Kuspel_1989}).}
\textit{Let a function $\Psi_{\beta}$ of the form \eqref{Psi_beta}, which generates classes $C_{\beta,p}^\psi$, $p=1,\infty$, satisfy the condition $C_{y,2n}$ for some $n\in\mathbb{N}$ when $y$ is a point at which the function $|(\Psi_\beta\ast\varphi_n)(t)|$, ${\varphi_n(t)=\textnormal{sgn}\sin nt}$, takes its maximum value.
Then }
\begin{equation*}
d_{2n}(C_{\beta,\infty}^\psi, C)\geq\|\Psi_\beta\ast\varphi_n\|_C,
\end{equation*}
\begin{equation*}
d_{2n-1}(C_{\beta,1}^\psi, L)\geq\|\Psi_\beta\ast\varphi_n\|_C.
\end{equation*}

Sufficient conditions of the inclusion $\Psi_\beta\in C_{y,2n}$ for kernels of the form \eqref{Psi_beta} were suggested in \cite{Kuspel_1988,Shevaldin_1992,Nguen_1994,Shevaldin_1993,Serdyuk_1998,Serdyuk_1995,Serdyuk_1999}.
By using these conditions, these authors succeeded to apply Theorem~1 and obtain sharp estimates for widths $d_m(C_{\beta,\infty}^\psi,C)$ and $d_m(C_{\beta,1}^\psi,L)$ in several new cases.

\section{Main results}
\label{}
Proceed to the presentation of main results of our paper.
For any fixed $q\in(0,1)$, let $n_q$ denote the smallest number $n\geq9$, for which the following inequality is satisfied:
\begin{equation}\label{umova_n_0}
\dfrac{43}{10(1-q)}q^{\sqrt{n}}+\dfrac{160}{57(n-\sqrt{n})}\; \dfrac{q}{(1-q)^2}\leq
\left(\dfrac{1}{2}+\dfrac{2q}{(1+q^2)(1-q)}\right)\left(\dfrac{1-q}{1+q}\right)^{\frac {4}{1-q^2}}.
\end{equation}

By means of the suggested notation, we can formulate the following statement.

\textbf{Theorem 2.}
\textit{Let $q\in(0,1)$.
Then equalities \eqref{dno1} and \eqref{dno2} are satisfied for any $\beta\in \mathbb{R}$ and all numbers $n\geq n_q$}.

\textbf{Proof.}
In accordance with Theorem~1, it is sufficient to prove that  Poisson kernels $P_{q,\beta}(t)$ satisfy the condition $C_{y_0,2n}$ for any $q\in(0,1)$, $\beta\in \mathbb{R}$, and all numbers $n\geq n_q$ (here we assume that $y_0$ is a point at which the function $|\Phi_{q,\beta,n}(\cdot)|$ takes its maximum value, $\Phi_{q,\beta,n}(\cdot)=(P_{q,\beta}\ast\varphi_n)(\cdot)$, and $\varphi_n(\cdot)$ is determined by equality~\eqref{varp}), i.e.,
\begin{equation*}
|\Phi_{q,\beta,n}(y_0)|=|(P_{q,\beta}\ast\varphi_n)(y_0)|=\|P_{q,\beta}\ast\varphi_n\|_C.
\end{equation*}

The function
\begin{equation*}
\Phi_{q,\beta,n}(\cdot)=(P_{q,\beta}\ast\varphi_n)(\cdot)=\dfrac{4}{\pi}\sum\limits_{\nu=0}^{\infty}\dfrac{q^{(2\nu+1)n}}{2\nu+1}\sin\left((2\nu+1)n\cdot-\dfrac{\beta\pi}{2}\right),
\end{equation*}
is periodic with period $2\pi/n$ and such that $\Phi_{q,\beta,n}(\cdot+\dfrac{\pi}{n})=-\Phi_{q,\beta,n}(\cdot)$.
Therefore, $\pi/n$-periodic function $|\Phi_{q,\beta,n}(\cdot)|$ on $[0,\dfrac{\pi}{n})$ takes its maximum value at a point $y_0=y_0(n,q,\beta)=\dfrac{\theta_n\pi}{n}$ where $\theta_n$ is the root of equation~\eqref{theta}, $\theta_n\in[0,1)$.
Since $|\Phi_{q,\beta,n}(\cdot)|=|\Phi_{q,\beta+2,n}(\cdot)|$, without loss of generality we can assume that $\beta\in[0,2)$.
For these values $\beta$, the unique root of equation \eqref{theta} on $[0,1)$ can be written down in the explicit way as follows:
\begin{equation}\label{theta_n}
\theta_n=1-[\beta]-\frac{1}{\pi}\arcsin\frac{(1-q^{2n})\cos\dfrac{\beta\pi}{2}}{\sqrt{1-2q^{2n}\cos\beta\pi+q^{4n}}},\;\;\beta\in[0,2),
\end{equation}
where $[a]$ is the integer part of $a$.

Let functions $\Psi_{\beta,1}(t)$ of the form \eqref{Ps}, which are generated by Poisson kernels ${\Psi_\beta(t)=P_{q,\beta}(t)}$, be denoted by $P_{q,\beta,1}(t)$, and the fundamental $SK$-spline $\overline{S\Psi}_{\beta,1}(y,\cdot)$ be denoted by $\overline{SP}_{q,\beta,1}(y,\cdot)$.
We use a representation of the function $(\overline{S\Psi}_{\beta,1}(y,t))_\beta^\psi$ from \cite{Serdyuk_1995} where it is proved that under the condition $|\lambda_j(y)|\not=0$, $j=\overline{1, n}$, for any $t\in(x_{k-1},x_k)$ we have the equality
\begin{equation}\label{SPsi_v0}
(\overline{S\Psi}_{\beta,1}(y,t))_\beta^\psi=\frac{\pi}{4n^2}\left(2\sum_{j=1}^{n-1}\dfrac{\sin jt_k\cdot\rho_j(y)-\cos jt_k\cdot\sigma_j(y)}{|\lambda_j(y)|^2\sin\dfrac{j\pi}{2n}}+\dfrac{(-1)^{k+1}\rho_n(y)}{|\lambda_n(y)|^2}\right),
\end{equation}
where
\begin{equation*}
\lambda_j(\cdot)=\dfrac1n \sum_{\nu=1}^{2n}e^{ij\nu\pi/n}\Psi_{\beta,1}(\cdot-\dfrac{\nu\pi}{n}),
\end{equation*}
$i$ is the imaginary unit, $\rho_j(\cdot)=\mathop{\text{Re}}(\lambda_j(\cdot))$, $\sigma_j(\cdot)=\mathop{\text{Im}}(\lambda_j(\cdot))$, $t_k=\dfrac{k \pi}{n}-\dfrac{\pi}{2n}$.

Changing the order of summation terms in the sum in the right-hand side of equality \eqref{SPsi_v0}, we obtain
\begin{equation*}
\sum_{j=1}^{n-1}\dfrac{\sin jt_k\cdot\rho_j(y)-\cos jt_k\cdot\sigma_j(y)}{|\lambda_j(y)|^2\sin\dfrac{j\pi}{2n}}=
\end{equation*}
\begin{equation*}
=\sum_{j=1}^{n-1}\dfrac{\sin (n-j)t_k\cdot\rho_{n-j}(y)-\cos (n-j)t_k\cdot\sigma_{n-j}(y)}{|\lambda_{n-j}(y)|^2\sin\dfrac{(n-j)\pi}{2n}}=
\end{equation*}
\begin{equation}\label{S_Psi}
=(-1)^{k+1}\sum_{j=1}^{n-1}\dfrac{\cos jt_k\cdot\rho_{n-j}(y)-\sin jt_k\cdot\sigma_{n-j}(y)}{|\lambda_{n-j}(y)|^2\cos\dfrac{j\pi}{2n}}.
\end{equation}

Let ${y=y_0}$. As it follows from \eqref{theta_n} (see also \cite[p.~538]{Stepanets_2005}), the following inclusions are valid: 
\begin{equation}\label{ny_0_1}
ny_0\in[\dfrac{\pi}{2},\pi) \text{ if } \beta\in[0,1)\cup[2,3),
\end{equation}
\begin{equation}\label{ny_0_2}
ny_0\in[0, \dfrac{\pi}{2}) \text{ if } {\beta\in[1,2)\cup[3,4)}.
\end{equation}
Then in accordance with equality (19) from the paper \cite{Stepanets_1994} and Lamma~2 from \cite{Stepanets_1994} we have $|\lambda_{j}(y_0)|>0$,
${j=\overline{1,n}}$. 
Thus, taking \eqref{SPsi_v0} and \eqref{S_Psi} into account, for fundamental $SK$-spline $\overline{S\Psi}_{\beta,1}(y,t)=\overline{SP}_{q,\beta,1}(y,t)$, which is generated by the Poisson kernel $P_{q,\beta}(t)$, we obtain the following representation:
\begin{equation*}
(\overline{SP}_{q,\beta,1}(y_0,t))_\beta^q=
\end{equation*}
\begin{equation}\label{SP_v0}
=\frac{(-1)^{k+1}\pi}{4n^2}\left(2\sum_{j=1}^{n-1}\dfrac{\cos jt_k\cdot\rho_{n-j}(y_0)-\sin jt_k\cdot\sigma_{n-j}(y_0)}{|\lambda_{n-j}(y_0)|^2\cos\dfrac{j\pi}{2n}}+\dfrac{\rho_n(y_0)}{|\lambda_n(y_0)|^2}\right),
\end{equation}
where
\begin{equation}\label{lambd}
\lambda_{l}(y_0)=\dfrac1n \sum_{\nu=1}^{2n}e^{il\nu\pi/n}P_{q,\beta,1}(y_0-\dfrac{\nu\pi}{n}),\;l=\overline{1,n}.
\end{equation}

By using equality \eqref{SP_v0}, we will write down a representation of the function $(\overline{SP}_{q,\beta,1}(y_0,t))_\beta^q$, which is more convenient for our further study.
To this end, let us prove first that values $\lambda_{n-j}(y_0)$ of type \eqref{lambd} for $j=\overline{0,n-1}$ can be expressed as follows:
\begin{equation}\label{lambda_n-j}
\lambda_{n-j}(y_0)=e^{-ijy_0}\left((-1)^s\left(\dfrac{q^{n-j}}{n-j}+\dfrac{q^{n+j}}{n+j}\right)+r_{j}(y_0)\right),
\end{equation}
where
\begin{align}
\label{r_n-j}
&r_{j}(y_0)=\sum_{\nu=1}^3 r_{j}^{(\nu)}(y_0),
\\
\nonumber
&r_{j}^{(1)}(y_0)=
\dfrac{q^{3n-j}e^{i(3ny_0-\frac{(\beta+1)\pi}{2})}}{3n-j}+
\\
\label{r_nj_1}
&+\sum\limits_{m=2}^{\infty}\left(\dfrac{q^{(2m+1)n-j}e^{i((2m+1)ny_0-\frac{(\beta+1)\pi}{2})}}{(2m+1)n-j}+
\dfrac{q^{(2m-1)n+j}e^{-i((2m-1)ny_0-\frac{(\beta+1)\pi}{2})}}{(2m-1)n+j}\right),
\\
\label{r_nj_2}
&r_{j}^{(2)}(y_0)=i\left(\dfrac{q^{n+j}}{n+j}-\dfrac{q^{n-j}}{n-j}\right)\cos(ny_0-\frac{\beta\pi}{2}),
\\
\label{r_nj_3}
&r_{j}^{(3)}(y_0)=(-1)^s\left(\dfrac{q^{n-j}}{n-j}+\dfrac{q^{n+j}}{n+j}\right)(|\sin(ny_0-\dfrac{\beta\pi}{2})|-1),
\end{align}
and the value $s=s(n,q,\beta)$ is determined by the equality
\begin{equation}\label{s}
(-1)^s=\mathop{\mathrm{sgn}} \sin(ny_0-\dfrac{\beta \pi}{2}).
\end{equation}

Note that equality \eqref{s} is true since due to \eqref{ny_0_1} and \eqref{ny_0_2} we have
\begin{equation*}
\sin(ny_0-\dfrac{\beta \pi}{2})=\sin ny_0\cos\dfrac{\beta \pi}{2}-\cos ny_0\sin\dfrac{\beta \pi}{2}\not=0.
\end{equation*}

Rewrite the kernel $P_{q,\beta,1}$ in the complex form
\begin{equation*}
P_{q,\beta,1}(t)=(P_{q,\beta}\ast B_1)(t)=\sum_{k=1}^\infty\frac{q^k}{k}\cos (kt-\frac{(\beta+1)\pi}{2})= \dfrac 12 \sideset{}{'}\sum\limits_{k=-\infty}^{\infty}c_ke^{ikt},
\end{equation*}
where
\begin{equation}\label{c_k}
c_k=\dfrac{q^k}{k}e^{-i\frac{(\beta+1)\pi}{2}},\;c_{-k}=\dfrac{q^k}{k}e^{i\frac{(\beta+1)\pi}{2}},\;k\in\mathbb{N},
\end{equation}
and the prime sign at the sum sign means that zero index summand in the summation is missed.

Substituting the kernel $P_{q,\beta,1}$ in \eqref{lambd} by the Fourier expansion of the kernel $P_{q,\beta,1}$, we obtain
\begin{equation*}
\lambda_l(y_0)=\dfrac1n \sum_{\nu=1}^{2n}e^{il\nu\pi/n}\dfrac 12 \sideset{}{'}\sum\limits_{k=-\infty}^{\infty}c_ke^{ik(y_0-\nu\pi/n)} =
\end{equation*}
\begin{equation*}
=\dfrac{1}{2n} \sum_{\nu=1}^{2n}\sideset{}{'}\sum\limits_{k=-\infty}^{\infty}c_k e^{i(ky_0+(l-k)\nu\pi/n)} =
\end{equation*}
\begin{equation}\label{eq:lambda_dop1}
=\dfrac{1}{2n} \sideset{}{'}\sum\limits_{k=-\infty}^{\infty}c_k e^{iky_0} \sum_{\nu=1}^{2n}e^{i((l-k)\nu\pi/n)}.
\end{equation}
It is easy to verify that 
\begin{equation}\label{eq:lambda_dop2}
\sum_{\nu=1}^{2n}e^{i((l-k)\nu\pi/n)}=\begin{cases}
0,& \text{if } k\not=l-2mn, m\in\mathbb{Z};\\
2n,& \text{if } k=l-2mn, m\in\mathbb{Z}.
\end{cases}
\end{equation}

Relations \eqref{eq:lambda_dop1} and \eqref{eq:lambda_dop2} for $l=\overline{1,n}$ imply the following representation:
\begin{equation*}
\lambda_l(y_0)=\sum\limits_{m=-\infty}^{+\infty}c_{l-2mn}e^{i(l-2mn)y_0}
=\sum\limits_{m=-\infty}^{+\infty}c_{2mn+l}e^{i(2mn+l)y_0}.
\end{equation*}
For $l=n-j\,$, $j=\overline{0,n-1}$, this implies that
\begin{equation*}
\lambda_{n-j}(y_0)=\sum\limits_{m=-\infty}^{+\infty}c_{(2m+1)n-j}e^{i((2m+1)n-j)y_0}=
\end{equation*}
\begin{equation}\label{lambda_n-j_y0}
=e^{-ijy_0}(c_{n-j}e^{iny_0}+c_{-(n+j)}e^{-iny_0}+r_{j}^{(1)}(y_0)).
\end{equation}
By using \eqref{c_k}, we can rewrite the first two summands in \eqref{lambda_n-j_y0} as follows:
\begin{equation*}
c_{n-j}e^{iny_0}+c_{-(n+j)}e^{-iny_0}=
\end{equation*}
\begin{equation*}
=\dfrac{q^{n-j}}{n-j}e^{i(ny_0-\frac{(\beta+1)\pi}{2})}+\dfrac{q^{n+j}}{n+j}e^{-i(ny_0-\frac{(\beta+1)\pi}{2})}=
\end{equation*}
\begin{equation*}
=\left(\dfrac{q^{n-j}}{n-j}+\dfrac{q^{n+j}}{n+j}\right)\cos(ny_0-\frac{(\beta+1)\pi}{2})+
\end{equation*}
\begin{equation*}
+i\left(\dfrac{q^{n-j}}{n-j}-\dfrac{q^{n+j}}{n+j}\right)\sin(ny_0-\frac{(\beta+1)\pi}{2})=
\end{equation*}
\begin{equation}\label{c_n-j-}
=\left(\dfrac{q^{n-j}}{n-j}+\dfrac{q^{n+j}}{n+j}\right)\sin(ny_0-\frac{\beta\pi}{2})+r_{j}^{(2)}(y_0).
\end{equation}
Representing $\sin(ny_0-\dfrac{\beta\pi}{2})$ in the form
\begin{equation*}\label{sin_ny}
\sin\left(ny_0-\dfrac{\beta\pi}{2}\right)=(-1)^{s}|\sin(ny_0-\dfrac{\beta\pi}{2})|,
\end{equation*}
from \eqref{c_n-j-} we obtain
\begin{equation*}
c_{n-j}e^{iny_0}+c_{-(n+j)}e^{-iny_0}=
(-1)^{s}\left(\dfrac{q^{n-j}}{n-j}+\dfrac{q^{n+j}}{n+j}\right)+
\end{equation*}
\begin{equation*}
+(-1)^{s}\left(\dfrac{q^{n-j}}{n-j}+\dfrac{q^{n+j}}{n+j}\right)(|\sin(ny_0-\dfrac{\beta\pi}{2})|-1)
+r_{j}^{(2)}(y_0)=
\end{equation*}
\begin{equation}\label{c_n-j}
=(-1)^{s}\left(\dfrac{q^{n-j}}{n-j}+\dfrac{q^{n+j}}{n+j}\right)+r_{j}^{(2)}(y_0)+r_{j}^{(3)}(y_0).
\end{equation}
Equalities \eqref{lambda_n-j_y0} and \eqref{c_n-j} imply formula \eqref{lambda_n-j}.

Using formulas \eqref{SP_v0} and \eqref{lambda_n-j}, we derive a new representation for the function $(\overline{SP}_{q,\beta,1}(y_0,t))_\beta^q$.

\textbf{Lemma 1.}
\textit{Let $q\in(0,1)$, $\beta\in \mathbb{R}$, $y_0=\dfrac{\theta_n\pi}{n}$, where $\theta_n$ is the unique root of equation \eqref{theta} on $[0,1)$.
Then for any $t\in(\dfrac{(k-1)\pi}{n},\dfrac{k\pi}{n})$, $k=\overline{1,2n}$, the following equality is satisfied}:
\begin{equation}\label{SP_Phi_}
(\overline{SP}_{q,\beta,1}(y_0,t))_\beta^q=
(-1)^{k+s+1}\frac{\pi}{4nq^{n}}\;
(\mathcal{P}_q(t_k-y_0)+\sum_{m=1}^5\gamma_m(y_0)),
\end{equation}
\textit{where} $\mathcal{P}_q(t)$ \textit{is the Poisson kernel for the heat conduction equation}
\begin{equation}\label{Phi}
\mathcal{P}_q(t)= \dfrac{1}{2}+2\sum_{j=1}^{\infty}\dfrac{\cos jt}{q^j+q^{-j}},
\end{equation}
\textit{and}
\begin{align}
\label{gamma_1}
\gamma_1(y_0)&=\gamma_1(k,y_0)=2\sum_{j=[\sqrt{n}]+1}^{n-1}\dfrac{\cos j(t_k-y_0)}{\dfrac{n}{q^{n}}|\lambda_{n-j}(y_0)|\cos\dfrac{j\pi}{2n}},
\\
\gamma_2(y_0)&=\gamma_2(k,y_0)= \nonumber
\\
\label{gamma_2}
&=(-1)^{s}\frac{q^{n}}{n}
\left(\dfrac{z_{0}(y_0)}{|\lambda_{n}(y_0)|^2}+
2\sum_{j=1}^{n-1}\dfrac{z_{j}(y_0)}{|\lambda_{n-j}(y_0)|^2\cos\dfrac{j\pi}{2n}}\right),
\\
\label{gamma_3}
\gamma_3(y_0)&=-\dfrac{R_0(y_0)\dfrac{n}{q^{n}}}{2(2+R_0(y_0)\dfrac{n}{q^{n}})},
\\
\label{gamma_4}
\gamma_4(y_0)&=\gamma_4(k,y_0)=-2\sum_{j=1}^{[\sqrt{n}]}\dfrac{\delta_{j}(y_0)\cos j(t_k-y_0)}{\dfrac{n}{q^{n}}|\lambda_{n-j}(y_0)|\cos\dfrac{j\pi}{2n}},
\\
\label{gamma_5}
\gamma_5(y_0)&=\gamma_5(k,y_0)=-2\sum\limits_{j=[\sqrt{n}]+1}^{\infty}\dfrac{\cos j(t_k-y_0)}{q^j+q^{-j}},
\\
\label{delta_0}
\delta_{j}(y_0)&=\dfrac{n|\lambda_{n-j}(y_0)|\cos\dfrac{j\pi}{2n}}{(q^{-j}+q^{j})q^n}-1,\;j=\overline{1,[\sqrt{n}]},
\\
z_{j}(y_0)&=|r_{j}(y_0)|\cos(j(t_k-y_0)+\arg(r_{j}(y_0)))+\nonumber
\\
\label{z_nj}
&+(-1)^{s+1}R_{j}(y_0)\cos(j(t_k-y_0)),\;j=\overline{0,n-1},
\\
\label{|lambda_n-j|}
R_{j}(y_0)&=|\lambda_{n-j}(y_0)|- \dfrac{q^{n-j}}{n-j}-\dfrac{q^{n+j}}{n+j},\;j=\overline{0,n-1},
\end{align}
\textit{$t_k=\dfrac{k \pi}{n}-\dfrac{\pi}{2n}$, and values $\lambda_{l}(y_0)$, $r_{j}(y_0)$, $j=\overline{0,n-1}$, and $s$ are determined by equalities \eqref{lambd}, \eqref{r_n-j}, and \eqref{s} respectively.}

\textbf{Proof.}
Transform the numerator of each term in the right-hand side of equality~\eqref{SP_v0}.
To this end, taking \eqref{lambda_n-j} into account, we write down the following equalities:
\begin{equation*}
\rho_{n-j}(y_0)= \mathop{\text{Re}}(\lambda_{n-j}(y_0))=
\end{equation*}
\begin{equation}\label{rho_n-j}
=(-1)^{s}\left(\dfrac{q^{n-j}}{n-j}+\dfrac{q^{n+j}}{n+j}\right)\cos jy_0+\mathop{\text{Re}}(e^{-ijy_0}r_{j}(y_0));
\end{equation}
\begin{equation*}
\sigma_{n-j}(y_0)=\mathop{\text{Im}}(\lambda_{n-j}(y_0))=
\end{equation*}
\begin{equation}\label{sigma_n-j}
=(-1)^{s+1}\left(\dfrac{q^{n-j}}{n-j}+\dfrac{q^{n+j}}{n+j}\right)\sin jy_0+\mathop{\text{Im}}(e^{-ijy_0}r_{j}(y_0)).
\end{equation}
From equalities \eqref{rho_n-j} and \eqref{sigma_n-j} we obtain
\begin{equation*}
\cos jt_k\cdot\rho_{n-j}(y_0)-\sin jt_k\cdot\sigma_{n-j}(y_0)=
\end{equation*}
\begin{equation*}
=\cos jt_k
\left((-1)^{s}\left(\dfrac{q^{n-j}}{n-j}+\dfrac{q^{n+j}}{n+j}\right)\cos jy_0+\mathop{\text{Re}}(e^{-ijy_0} r_{j}(y_0))\right)+
\end{equation*}
\begin{equation*}
+\sin jt_k
\left((-1)^{s}\left(\dfrac{q^{n-j}}{n-j}+\dfrac{q^{n+j}}{n+j}\right)\sin jy_0-\mathop{\text{Im}}(e^{-ijy_0} r_{j}(y_0))\right)=
\end{equation*}
\begin{equation*}
=
(-1)^{s}\left(\dfrac{q^{n-j}}{n-j}+\dfrac{q^{n+j}}{n+j}\right)\cos(j(t_k-y_0))+
\end{equation*}
\begin{equation*}
+\cos jt_k\cdot\mathop{\text{Re}} (e^{-ijy_0} r_{j}(y_0))-\sin jt_k\cdot\mathop{\text{Im}} (e^{-ijy_0} r_{j}(y_0))=
\end{equation*}
\begin{equation*}
=
(-1)^{s}\left(\dfrac{q^{n-j}}{n-j}+\dfrac{q^{n+j}}{n+j}+R_{j}(y_0)\right)\cos(j(t_k-y_0))+z_{j}(y_0)=
\end{equation*}
\begin{equation}\label{cos+sin}
=
(-1)^{s}|\lambda_{n-j}(y_0)|\cos(j(t_k-y_0))+z_{j}(y_0),
\end{equation}
where
\begin{equation*}
z_{j}(y_0)=\cos jt_k\cdot\mathop{\text{Re}}(e^{-ijy_0} r_{j}(y_0))-\sin jt_k\cdot\mathop{\text{Im}} (e^{-ijy_0} r_{j}(y_0))+
\end{equation*}
\begin{equation*}
+(-1)^{s+1}R_{j}(y_0)\cos(j(t_k-y_0)),
\end{equation*}
and $R_{j}(y_0)$ are defined in \eqref{|lambda_n-j|}.
Due to the obvious equality
\begin{equation*}
e^{-ijy_0} r_{j}(y_0)=|r_{j}(y_0)|(\cos (\arg(r_{j}(y_0))-jy_0)+i\sin (\arg(r_{j}(y_0))-jy_0))
\end{equation*}
we can represent the quantity $z_{j}(y_0)$ in the form of \eqref{z_nj}.

For $j=0$, formula \eqref{lambda_n-j} turns into the equality
\begin{equation}\label{lambda_n}
\lambda_{n}(y_0)=(-1)^{s}\, 2\dfrac{q^{n}}{n}+r_{0}(y_0),
\end{equation}
where $r_{0}(y_0)$ is determined by formula \eqref{r_n-j}, and in that formula we have
\begin{align}
\label{r_n1}
&r_{0}^{(1)}(y_0)=2\sum\limits_{m=2}^{\infty}\dfrac{q^{(2m-1)n}}{(2m-1)n}\cos((2m-1)ny_0-\frac{(\beta+1)\pi}{2}),
\\
\label{r_n2}
&r_{0}^{(2)}(y_0)=0,
\\
\label{r_n3}
&r_{0}^{(3)}(y_0)=(-1)^{s}\,2\dfrac{q^{n}}{n} (|\sin(ny_0-\dfrac{\beta\pi}{2})|-1).
\end{align}

Relations \eqref{lambda_n}--\eqref{r_n3} imply that $\sigma_{n}(y_0)=0$ and hence
\begin{equation*}
\rho_{n}(y_0)=\lambda_{n}(y_0)=(-1)^{s}\,2\dfrac{q^{n}}{n}+r_{0}(y_0).
\end{equation*}
From this, by using \eqref{z_nj} and \eqref{|lambda_n-j|}, we can conclude that
\begin{equation*}
\rho_{n}(y_0)=(-1)^{s}\left(2\dfrac{q^{n}}{n}+R_0(y_0)\right)+z_{0}(y_0)=
\end{equation*}
\begin{equation}\label{rho_n}
=(-1)^{s}|\lambda_{n}(y_0)|+z_{0}(y_0),
\end{equation}
where
\begin{equation*}
z_0(y_0)=r_0(y_0)+(-1)^{s+1} R_0(y_0).
\end{equation*}
By using representation \eqref{SP_v0} and equalities \eqref{cos+sin} and \eqref{rho_n}, we obtain
\begin{equation*}
(\overline{SP}_{q,\beta,1}(y_0,t))_\beta^q=
\end{equation*}
\begin{equation*}
=\frac{(-1)^{k+1}\pi}{4nq^n}
\left((-1)^{s}
\left(
2\frac{q^n}{n}\sum_{j=1}^{[\sqrt{n}]}\dfrac{\cos j(t_k-y_0)}{|\lambda_{n-j}(y_0)|\cos\dfrac{j\pi}{2n}}
+\dfrac{q^n}{n|\lambda_n(y_0)|}+
\right.\right.
\end{equation*}
\begin{equation*}
\left.\left.
+2\frac{q^n}{n}\sum_{j=[\sqrt{n}]+1}^{n-1}\dfrac{\cos j(t_k-y_0)}{|\lambda_{n-j}(y_0)|\cos\dfrac{j\pi}{2n}}\right)
+2\frac{q^n}{n}\sum_{j=1}^{n-1}\dfrac{z_{j}(y_0)}{|\lambda_{n-j}(y_0)|^2\cos\dfrac{j\pi}{2n}}
+\dfrac{q^nz_{0}(y_0)}{n|\lambda_{n}(y_0)|^2}\right)=
\end{equation*}
\begin{equation*}
=\frac{(-1)^{k+s+1}\pi}{4nq^n}
\times
\end{equation*}
\begin{equation}\label{Phi_x}
\times
\left(2\frac{q^n}{n}\sum_{j=1}^{[\sqrt{n}]}\dfrac{\cos j(t_k-y_0)}{|\lambda_{n-j}(y_0)|\cos\dfrac{j\pi}{2n}}
+\dfrac{q^n}{n|\lambda_n(y_0)|}+\gamma_1(y_0)+\gamma_2(y_0)
\right).
\end{equation}
Due to \eqref{|lambda_n-j|}, we have
\begin{equation}\label{s0}
\dfrac{q^{n}}{n|\lambda_n(y_0)|}=\dfrac{1}{2+R_0(y_0)\dfrac{n}{q^{n}}}=\dfrac{1}{2}-\dfrac{R_0(y_0)\dfrac{n}{q^{n}}}{2(2+R_0(y_0)\dfrac{n}{q^{n}})}=\dfrac{1}{2}+\gamma_3(y_0).
\end{equation}
Taking formula \eqref{delta_0} into account, we can write down the following equalities:
\begin{equation*}
2\dfrac{q^{n}}{n}\sum_{j=1}^{[\sqrt{n}]}\dfrac{\cos j(t_k-y_0)}{|\lambda_{n-j}(y_0)|\cos\dfrac{j\pi}{2n}}=
2\sum_{j=1}^{[\sqrt{n}]}\dfrac{\cos j(t_k-y_0)}{(q^j+q^{-j})(1+\delta_{j}(y_0))}=
\end{equation*}
\begin{equation*}
=2\sum_{j=1}^{[\sqrt{n}]}\dfrac{\cos j(t_k-y_0)}{q^j+q^{-j}}-2\sum_{j=1}^{[\sqrt{n}]}\dfrac{\delta_{j}(y_0)\cos j(t_k-y_0)}{(q^j+q^{-j})(1+\delta_{j}(y_0))}=
\end{equation*}
\begin{equation}\label{s1}
=\mathcal{P}_q(t_k-y_0)-\frac{1}{2}+\gamma_4(y_0)+\gamma_5(y_0),
\end{equation}
Relations \eqref{Phi_x}--\eqref{s1} imply \eqref{SP_Phi_}.
The lemma is proved.

The following lemma contains a lower estimate for the minimum value of the kernel $\mathcal{P}_q(\cdot)$ of the form~\eqref{Phi}.

\textbf{Lemma 2.}
\textit{Let $q\in(0,1)$.
Then for any $x\in\mathbb{R}$ the following inequality is satisfied}:
\begin{equation}\label{f_x}
\mathcal{P}_q(x)>\left(\dfrac{1}{2}+\dfrac{2q}{(1+q^2)(1-q)}\right)\left(\dfrac{1-q}{1+q}\right)^{\frac {4}{1-q^2}}.
\end{equation}

\textbf{Proof.}
We use the following representation of the kernel $\mathcal{P}_q(x)$, which is derived in the theory of elliptic functions (see e.g. \cite[p.~867]{Gradshteyn_2007}):
\begin{equation}\label{qr1}
\mathcal{P}_q(x)=\frac{K}{\pi}\mathop{\text{dn}}\left(\frac{Kx}{\pi}\right),
\end{equation}
where
\begin{equation}\label{K_}
K=\pi\left(\dfrac{1}{2}+2\sum_{j=1}^{\infty}\dfrac{q^j}{1+q^{2j}}\right).
\end{equation}

In accordance with formula (8.146.22) from \cite[p.~868]{Gradshteyn_2007}, we have
\begin{equation}\label{qr2}
\mathop{\text{dn}}\left(\frac{Kx}{\pi}\right)=\exp\left\{ -8\sum\limits_{j=1}^{\infty}\dfrac{1}{2j-1}\,\dfrac{q^{2j-1}}{1-q^{2(2j-1)}}\sin^2\dfrac {2j-1}{2}x\right\},
\end{equation}
Due to the equality
\begin{equation*}
\sum\limits_{\nu=1}^{\infty}\dfrac{q^{2\nu-1}}{2\nu-1}=\dfrac 12\ln\dfrac{1+q}{1-q}
\end{equation*}
(see e.g. \cite[p.~53]{Gradshteyn_2007}), the following estimate is valid:
\begin{equation*}
\sum\limits_{j=1}^{\infty}\dfrac{1}{2j-1}\,\dfrac{q^{2j-1}}{1-q^{2(2j-1)}}\sin^2\dfrac {2j-1}{2}x<
\end{equation*}
\begin{equation}\label{sum}
<\dfrac {1}{1-q^2}\sum\limits_{j=1}^{\infty}\dfrac{q^{2j-1}}{2j-1}
=\dfrac {1}{1-q^2}\,\dfrac 12\ln\dfrac{1+q}{1-q}.
\end{equation}
Relations \eqref{qr2} and \eqref{sum} imply the inequality
\begin{equation}\label{dn}
\mathop{\text{dn}}\left(\frac{Kx}{\pi}\right)>\left(\dfrac{1-q}{1+q}\right)^{\frac {4}{1-q^2}},
\end{equation}
and \eqref{K_} implies the inequality
\begin{equation}\label{K}
\frac{K}{\pi}> \dfrac{1}{2}+\dfrac{2}{1+q^2}\sum_{j=1}^{\infty}q^j= \dfrac{1}{2}+\dfrac{2q}{(1+q^2)(1-q)}.
\end{equation}
Formulas \eqref{qr1}, \eqref{dn}, and \eqref{K} imply \eqref{f_x}.
The lemma is proved.

The following statement contains an upper estimate for the sum $\sum\limits_{k=1}^5|\gamma_k(y_0)|$.

\textbf{Lemma 3.}
\textit{Let $q\in(0,1)$, $\beta\in\mathbb{R}$, and values $\gamma_k(y_0)$, $k=\overline{1,5}$, be determined by equalities \eqref{gamma_1}--\eqref{gamma_5}.
Then for $n\geq9$ under the condition}
\begin{equation}\label{umova_z}
\dfrac{q^{n}}{1-q^{2n}}\leq \dfrac {7q^{\sqrt{n}}}{37n^2}
\end{equation}
\textit{the following estimate is valid}:
\begin{equation*}
\sum\limits_{k=1}^5|\gamma_k(y_0)|\leq\dfrac{43}{10(1-q)}q^{\sqrt{n}}+\dfrac{160}{57(n-\sqrt{n})}\; \dfrac{q}{(1-q)^2}.
\end{equation*}

\textbf{Proof.}
Let us estimate each summand $|\gamma_k(y_0)|$, $k=\overline{1,5}$, separately.
To this end, at first, we find upper estimates for values $|r_{j}(y_0)|$ and $|R_{j}(y_0)|$ with $j=\overline{0,n-1}$.
Given that due to convexity of the sequence $\dfrac{q^k}{k}$ we have the inequality ${\dfrac{q^{k-j}}{k-j}+\dfrac{q^{k+j}}{k+j}<\dfrac{q^{k-n}}{k-n}+\dfrac{q^{k+n}}{k+n}}$, $k>n$, $j=\overline{0,n-1}$, from \eqref{r_nj_1} we obtain
\begin{equation*}
|r_{j}^{(1)}(y_0)|\leq
\dfrac{q^{3n-j}}{3n-j}+\sum\limits_{m=2}^{\infty}\left(\dfrac{q^{(2m+1)n-j}}{(2m+1)n-j}+\dfrac{q^{(2m-1)n+j}}{(2m-1)n+j}\right)=
\end{equation*}
\begin{equation*}
=\sum\limits_{m=1}^{\infty}\left(\dfrac{q^{(2m+1)n-j}}{(2m+1)n-j}+\dfrac{q^{(2m+1)n+j}}{(2m+1)n+j}\right)\leq
\end{equation*}
\begin{equation*}
\leq\sum\limits_{m=1}^{\infty}\left(\dfrac{q^{2mn}}{2mn}+\dfrac{q^{2(m+1)n}}{2(m+1)n}\right)=
\end{equation*}
\begin{equation}\label{|r_nj1|}
= \dfrac{q^{2n}}{2n}+\sum\limits_{m=2}^{\infty}\dfrac{q^{2mn}}{mn}\leq
\dfrac{1}{2n}\sum\limits_{m=1}^{\infty}q^{2mn}=\dfrac{q^{2n}}{2n(1-q^{2n})}.
\end{equation}
Since $y_0=\dfrac{\theta_n\pi}{n}$, from \eqref{theta} we derive
\begin{equation*}
\left|\cos\left(ny_0-\dfrac{\beta\pi}{2}\right)\right|=\left|\sum\limits_{\nu=1}^{\infty}q^{2\nu n}\cos\left((2\nu+1)ny_0-\dfrac{\beta\pi}{2}\right)\right|\leq
\end{equation*}
\begin{equation}\label{cos_y}
\leq \sum\limits_{\nu=1}^{\infty}q^{2\nu n}=\dfrac{q^{2n}}{1-q^{2n}}.
\end{equation}
From \eqref{r_nj_2} and \eqref{cos_y} we have
\begin{equation}\label{|r_nj2|}
|r_{j}^{(2)}(y_0)|\leq|\cos(ny_0-\frac{\beta\pi}{2})|\left(\dfrac{q^{n-j}}{n-j}-\dfrac{q^{n+j}}{n+j}\right)
\leq\dfrac{q^{2n}}{1-q^{2n}}\left(q-\dfrac{q^{2n-1}}{2n-1}\right).
\end{equation}
Relation \eqref{cos_y} implies that
\begin{equation}\label{a_n}
0\leq 1-|\sin(ny_0-\dfrac{\beta\pi}{2})|\leq |\cos(ny_0-\dfrac{\beta\pi}{2})|\leq \dfrac{q^{2n}}{1-q^{2n}}.
\end{equation}
From \eqref{r_nj_3} and \eqref{a_n} we conclude that
\begin{equation}\label{|alpha_n|}
|r_{j}^{(3)}(y_0)|\leq\dfrac{q^{2n}}{1-q^{2n}}\left(q+\dfrac{q^{2n-1}}{2n-1}\right).
\end{equation}
Hence, for the value $r_{j}(y_0)$, from \eqref{|r_nj1|}, \eqref{|r_nj2|}, and \eqref{|alpha_n|} we obtain the estimate
\begin{equation}\label{|r_nj|}
|r_{j}(y_0)|\leq |\sum\limits_{\nu=1}^{3}r_{j}^{(\nu)}(y_0)|
 \leq \dfrac{q^{2n}}{1-q^{2n}}\left(2q+\dfrac{1}{2n}\right)\leq \dfrac{37q^{2n}}{18(1-q^{2n})},j=\overline{0,n-1}.
\end{equation}

For $j=0$, estimate \eqref{|r_nj|} can be improved.
Indeed, due to \eqref{r_n1} and \eqref{r_n3} for $j=0$, we have
\begin{equation*}
|r_{0}^{(1)}(y_0)|\leq2\sum\limits_{m=2}^{\infty}\dfrac{q^{(2m-1)n}}{(2m-1)n}\leq
\dfrac{2}{3n}\sum\limits_{m=2}^{\infty}q^{(2m-1)n}=\dfrac{2}{3n}\dfrac{q^{3n}}{1-q^{2n}},
\end{equation*}
\begin{equation*}
\left| r_{0}^{(3)}(y_0)  \right|\leq \dfrac{2q^{3n}}{n(1-q^{2n})}.
\end{equation*}
Then, taking \eqref{r_n2} into account, we obtain
\begin{equation}\label{|r_n|}
|r_{0}(y_0)|\leq |r_{0}^{(1)}(y_0)+r_{0}^{(3)}(y_0)| \leq \dfrac{8}{3n} \dfrac{q^{3n}}{1-q^{2n}}.
\end{equation}

For the value $|R_{j}(y_0)|$ of the form \eqref{|lambda_n-j|}, relation \eqref{lambda_n-j} implies the representation
\begin{equation*}
|\lambda_{n-j}(y_0)|=\left|(-1)^{s}\left(\dfrac{q^{n-j}}{n-j}+\dfrac{q^{n+j}}{n+j}\right)+r_{j}(y_0)\right|,
\end{equation*}
which, in its turn, immediately implies the following estimates:
\begin{equation}\label{|lambda_n-j|0'}
|\lambda_{n-j}(y_0)|\leq \dfrac{q^{n-j}}{n-j}+\dfrac{q^{n+j}}{n+j}+|r_{j}(y_0)|,
\end{equation}
\begin{equation}\label{|lambda_n-j|0''}
|\lambda_{n-j}(y_0)|\geq \dfrac{q^{n-j}}{n-j}+\dfrac{q^{n+j}}{n+j}-|r_{j}(y_0)|.
\end{equation}
Due to \eqref{|lambda_n-j|}, \eqref{|lambda_n-j|0'}, and \eqref{|lambda_n-j|0''}, we have
\begin{equation}\label{|R_n,j|}
|R_{j}(y_0)|\leq |r_{j}(y_0)|,\; j=\overline{0,n-1}.
\end{equation}

Now let us estimate the value $|\gamma_1(y_0)|$.
Since for $x\in[0,\dfrac \pi2)$ the inequality $\cos x\geq 1-\dfrac {2x}{\pi}>0$ holds, we obtain
\begin{equation}\label{cos1}
\cos \dfrac{j\pi}{2n}\geq 1-\dfrac {j}{n}=\dfrac {n-j}{n},\; j=\overline{0,n-1}.
\end{equation}
Taking estimates \eqref{|lambda_n-j|0''}, \eqref{|r_nj|}, and \eqref{cos1} into account, we conclude that
\begin{equation*}
\dfrac{n}{q^{n}}|\lambda_{n-j}(y_0)|\cos\dfrac{j\pi}{2n}\geq \dfrac{n}{q^{n}} \left(\dfrac{q^{n-j}}{n-j}+\dfrac{q^{n+j}}{n+j}-\dfrac{37q^{2n}}{18(1-q^{2n})}\right)\dfrac {n-j}{n}=
\end{equation*}
\begin{equation}\label{|lambda_n-j|3}
=q^{-j}+\dfrac{n-j}{n+j}q^{j}-\dfrac{37(n-j)q^{n}}{18(1-q^{2n})}.
\end{equation}
Since for $n\geq9$ we have $\dfrac{9q^{-j}}{10}>\dfrac{9}{10}>\dfrac{7}{9}>\dfrac{7q^{\sqrt{n}}}{n}$, condition \eqref{umova_z} implies the inequality
\begin{equation*}
\dfrac{9q^{-j}}{370n}>\dfrac{q^{n}}{1-q^{2n}},
\end{equation*}
and this inequality is equivalent to the following one:
\begin{equation}\label{um1}
\dfrac{q^{-j}}{20}>\dfrac{37nq^{n}}{18(1-q^{2n})},\; j=\overline{0,n-1}.
\end{equation}
Due to \eqref{um1} the following estimates are satisfied:
\begin{equation*}
q^{-j}+\dfrac{n-j}{n+j}q^{j}-\dfrac{37(n-j)q^{n}}{18(1-q^{2n})}=
\end{equation*}
\begin{equation}\label{ots}
=\frac{19}{20}q^{-j}+\frac{q^{-j}}{20}+\dfrac{n-j}{n+j}q^{j}-\dfrac{37(n-j)q^{n}}{18(1-q^{2n})}>\frac{19}{20}q^{-j},\; j=\overline{0,n-1}.
\end{equation}
Combining \eqref{|lambda_n-j|3} and \eqref{ots}, we see that
\begin{equation}\label{|lambda_n-j|2}
\dfrac{n}{q^{n}}|\lambda_{n-j}(y_0)|\cos\dfrac{j\pi}{2n}\geq \dfrac{19q^{-j}}{20}.
\end{equation}
Therefore, by using \eqref{|lambda_n-j|2}, from \eqref{gamma_1} we obtain
\begin{equation}\label{r3}
\left|\gamma_1(y_0)\right|\leq
\dfrac{40}{19}\sum_{j=[\sqrt{n}]+1}^{n-1}q^j=\dfrac{40(q^{[\sqrt{n}]+1}-q^{n-1})}{19(1-q)}\leq
\dfrac{40q^{\sqrt{n}}}{19(1-q)}.
\end{equation}

Let us estimate the value $|\gamma_2(y_0)|$.
Taking estimates \eqref{|r_nj|}, \eqref{|lambda_n-j|0''}, \eqref{cos1}, and \eqref{ots} into account, we have
\begin{equation*}
\dfrac{n}{q^{n}}|\lambda_{n-j}(y_0)|^2\cos\dfrac{j\pi}{2n}\geq \dfrac{n}{q^{n}} \left(\dfrac{q^{n-j}}{n-j}+\dfrac{q^{n+j}}{n+j}-\dfrac{37q^{2n}}{18(1-q^{2n})}\right)^2\dfrac {n-j}{n}=
\end{equation*}
\begin{equation}\label{|lambda_n-j|_2}
=\dfrac{q^{n}}{n-j} \left(q^{-j}+\dfrac{n-j}{n+j}q^{j}-\dfrac{37(n-j)q^{n}}{18(1-q^{2n})}\right)^2>
\dfrac{361q^{n-2j}}{400n}.
\end{equation}
Relations \eqref{|R_n,j|} and \eqref{z_nj} imply $|z_{j}(y_0)|\leq 2|r_{j}(y_0)|$.
Therefore, due to \eqref{|r_nj|}, \eqref{|lambda_n-j|_2}, and condition \eqref{umova_z}, from \eqref{gamma_2} we obtain
\begin{equation*}
\left|\gamma_2(y_0)\right|\leq
\dfrac{1600}{361}\max_{0\leq j\leq n-1}| r_{j}(y_0)| \dfrac{n}{q^{n}}\sum_{j=0}^{n-1}q^{2j}<
\dfrac{29600\,nq^{n}}{3249(1-q^{2n})}\sum_{j=0}^{\infty}q^{2j}\leq
\end{equation*}
\begin{equation}\label{r1}
\leq\dfrac{29600\,n}{3249}\;\dfrac{7q^{\sqrt{n}}}{37n^2}\;\dfrac{1}{1-q^2}  \leq
\dfrac {5600}{3249n} q^{\sqrt{n}}\;\dfrac{1}{1-q^2} .
\end{equation}

Let us estimate the value $|\gamma_3(y_0)|$.
Condition \eqref{umova_z} for $n\geq9$ implies the inequality
\begin{equation*}
q^{2n}\leq \dfrac{49}{8982009}.
\end{equation*}
Then by using \eqref{gamma_3}, \eqref{|R_n,j|}, \eqref{|r_n|} and \eqref{umova_z}, we obtain
\begin{equation*}
|\gamma_3(y_0)|
\leq\dfrac{ \dfrac{8q^{2n}}{3(1-q^{2n})}}{2\left|2-\dfrac{8q^{2n}}{3(1-q^{2n})}\right|}=\dfrac{2q^{2n}}{3-7q^{2n}}
=\dfrac{1-q^{2n}}{3-7q^{2n}}\;\dfrac{2q^{2n}}{1-q^{2n}}=
\end{equation*}
\begin{equation*}
=\left(\dfrac 17+\dfrac{4}{7(3-7q^{2n})}\right)\dfrac{2q^{2n}}{1-q^{2n}}<
\end{equation*}
\begin{equation}\label{r2_2}
<\dfrac 37 \dfrac{2q^{2n}}{1-q^{2n}}<\dfrac{6q^{n+\sqrt{n}}}{37n^2}.
\end{equation}

Before we estimate the value $|\gamma_4(y_0)|$, we will obtain upper estimates for $|\delta_j(y_0)|$ of the form \eqref{delta_0}.

Due to \eqref{|lambda_n-j|} we have
\begin{equation*}
\dfrac{n}{q^{n}}|\lambda_{n-j}(y_0)|\cos\dfrac{j\pi}{2n}=
\end{equation*}
\begin{equation*}
=\left(\dfrac{n}{n-j}\dfrac{q^{n-j}}{q^{n}}+\dfrac{n}{n+j}\dfrac{q^{n+j}}{q^{n}}+ R_{j}(y_0)\dfrac{n}{q^{n}}\right)\cos\dfrac{j\pi}{2n}=
\end{equation*}
\begin{equation*}
=\left((1+\dfrac{j}{n-j})q^{-j}+(1-\dfrac{j}{n+j})q^{j}+R_{j}(y_0)\dfrac{n}{q^{n}}\right)\cos\dfrac{j\pi}{2n}=
\end{equation*}
\begin{equation*}
=(q^{-j}+q^{j})(1-2\sin^2\dfrac{j\pi}{4n})+
\end{equation*}
\begin{equation}\label{cos_lambda0}
+\left(\dfrac{j}{n-j}q^{-j}-\dfrac{j}{n+j}q^{j}+R_{j}(y_0)\dfrac{n}{q^{n}}\right)\cos\dfrac{j\pi}{2n}.
\end{equation}

Due to relations \eqref{delta_0}, \eqref{|r_nj|}, \eqref{|R_n,j|}, \eqref{cos_lambda0}, and the fact that the sequence $q^k$ is convex, for values $|\delta_{j}(y_0)|$ we obtain the following inequalities:
\begin{equation*}
|\delta_{j}(y_0)|\leq 2\sin^2\dfrac{j\pi}{4n}
+\dfrac{1}{q^{-j}+q^{j}}\left(\dfrac{j}{n-j}q^{-j}+\dfrac{j}{n-j}q^{j}+|R_{j}(y_0)|\dfrac{n}{q^{n}}\right) \leq
\end{equation*}
\begin{equation*}
\leq 2\left(\dfrac{j\pi}{4n}\right)^2
+ \dfrac{j}{n-j}+\dfrac{n|r_{j}(y_0)|}{q^{n-j}+q^{n+j}}
\leq \dfrac{j^2\pi^2}{8n^2}
+ \dfrac{j}{n-j}+\dfrac{37nq^{n}}{36(1-q^{2n})}=
\end{equation*}
\begin{equation}\label{delta0}
= \dfrac{4j}{3(n-j)}+\left(\dfrac{j^2\pi^2}{8n^2}
+\dfrac{37nq^{n}}{36(1-q^{2n})}-\dfrac{j}{3(n-j)}\right).
\end{equation}
Let us prove that for all $j=\overline{1,[\sqrt{n}]}$ the expression in parentheses in the right-hand side of \eqref{delta0} is negative, i.e.,
\begin{equation}\label{umova4}
\dfrac{j}{3(n-j)}-\dfrac{j^2\pi^2}{8n^2}>\dfrac{37nq^{n}}{36(1-q^{2n})}.
\end{equation}
Then \eqref{delta0} and \eqref{umova4} will imply
\begin{equation}\label{delta}
|\delta_{j}(y_0)|\leq \dfrac{4j}{3(n-j)}.
\end{equation}

Indeed, for any fixed ${x\geq9}$, the function $f(x,\tau)=\dfrac{\tau}{3(x-\tau)}-\dfrac{\tau^2\pi^2}{8x^2}$ is convex upwards on $[1,\sqrt{x}]$.
Therefore, it takes its minimum value either at $\tau=1$ or at $\tau=\sqrt{x}$.
Consider the difference $f(x,1)-f(x,\sqrt{x})$ for $x\geq9$:
\begin{equation*}
f(x,1)-f(x,\sqrt{x})=\dfrac{1}{3(x-1)}-\dfrac{\pi^2}{8x^2}-\dfrac{\sqrt{x}}{3(x-\sqrt{x})}+\dfrac{x\pi^2}{8x^2}=
\end{equation*}
\begin{equation}\label{f-f}
=\dfrac{-8x^2\sqrt{x}+3\pi^2(x-1)^2}{24x^2(x-1)}.
\end{equation}
It is easy to verify that the function $g(x)=-8x^2\sqrt{x}+3\pi^2(x-1)^2$ is decreasing on $[9,+\infty)$ and $g(9)<0$. Therefore, due to \eqref{f-f} we have $f(x,1)-f(x,\sqrt{x})<0$.
For $n\geq9$, taking \eqref{umova_z} into account, we obtain for all $j=\overline{1,[\sqrt{n}]}$
\begin{equation*}
\dfrac{j}{3(n-j)}-\dfrac{j^2\pi^2}{8n^2}\geq
\dfrac{1}{3(n-1)}-\dfrac{\pi^2}{8n^2}>\dfrac{1}{3n}\left(1-\dfrac{\pi^2}{24}\right)>
\end{equation*}
\begin{equation*}
>\dfrac{7}{36n}>\dfrac{7q^{\sqrt{n}}}{36n}>\dfrac{37nq^{n}}{36(1-q^{2n})}.
\end{equation*}
This implies \eqref{umova4}, and hence, \eqref{delta} is also implied.

Formulas \eqref{gamma_4}, \eqref{|lambda_n-j|2}, and \eqref{delta} imply the following estimate of the value $\gamma_4(y_0)$ for $n\geq9$:
\begin{equation*}
|\gamma_4(y_0)|\leq
2\sum_{j=1}^{[\sqrt{n}]}\dfrac{\dfrac{4j}{3(n-j)}}{\dfrac{19q^{-j}}{20}}
=\frac{160}{57}\sum_{j=1}^{[\sqrt{n}]}\dfrac{j}{n-j}q^j\leq
\end{equation*}
\begin{equation*}
\leq
\dfrac{160}{57(n-\sqrt{n})}\sum_{j=1}^{[\sqrt{n}]}jq^j<\dfrac{160}{57(n-\sqrt{n})}\sum_{j=1}^{\infty}jq^j<
\end{equation*}
\begin{equation}\label{|R|}
<\dfrac{160}{57(n-\sqrt{n})}\; \dfrac{q}{(1-q)^2}.
\end{equation}

Due to \eqref{gamma_5}, we have the following estimate for the value $|\gamma_5(y_0)|$:
\begin{equation}\label{r4}
\left|\gamma_5(y_0)\right|\leq 2\sum\limits_{j=[\sqrt{n}]+1}^{\infty}q^j=2\dfrac{q^{[\sqrt{n}]+1}}{1-q}< 2\dfrac{q^{\sqrt{n}}}{1-q}.
\end{equation}

Taking estimates \eqref{r3}, \eqref{r1}, \eqref{r2_2}, \eqref{|R|}, and \eqref{r4} into account, for $n\geq9$ we obtain
\begin{equation*}
\sum_{k=1}^5|\gamma_k(y_0)|<
\end{equation*}
\begin{equation*}
<\dfrac{40q^{\sqrt{n}}}{19(1-q)}+
\dfrac {5600}{3249n} q^{\sqrt{n}}\;\dfrac{1}{1-q^2}+\dfrac{6q^{n+\sqrt{n}}}{37n^2}+\dfrac{160}{57(n-\sqrt{n})}\; \dfrac{q}{(1-q)^2}+\dfrac{2q^{\sqrt{n}}}{1-q}
<
\end{equation*}
\begin{equation*}
<\dfrac{q^{\sqrt{n}}}{1-q}(2.1053+0.1916+0.0021+2)+\dfrac{160}{57(n-\sqrt{n})}\; \dfrac{q}{(1-q)^2}<
\end{equation*}
\begin{equation*}
<\dfrac{43}{10(1-q)}q^{\sqrt{n}}+\dfrac{160}{57(n-\sqrt{n})}\; \dfrac{q}{(1-q)^2}.
\end{equation*}
The lemma is proved.

With proven above lemmas 2--3, under conditions \eqref{umova_n_0} and \eqref{umova_z} for $n\geq9$ we have
\begin{equation}\label{geq0}
{\mathcal{P}_q(t_k-y_0)+\sum\limits_{m=1}^5\gamma_m(y_0)}\geq0.
\end{equation}
Due to representation \eqref{SP_Phi_} and inequality \eqref{geq0}, we can conclude that under conditions \eqref{umova_n_0} and \eqref{umova_z} for $n\geq9$ the inclusion $P_{q,\beta}(t)\in C_{y_0,2n}$ is valid.

Note that for $q\in(0,\dfrac{9}{25}]$ condition \eqref{umova_z} is satisfied for any ${n\geq9}$.
To verify this, it suffices to observe that the sequence $\xi(n)=(n-\sqrt{n})\ln\dfrac{9}{25}+2\ln n-\ln\left(\dfrac{7}{37}\bigg(1-\left(\dfrac{3}{5}\right)^{36}\bigg)\right)$ is monotonously decreasing for $n\geq9$ and $\xi(9)<0$.
Therefore, for $n\geq9$ we have
\begin{equation}\label{q003}
(n-\sqrt{n})\ln\dfrac{9}{25}+2\ln n-\ln\left(\dfrac{7}{37}\bigg(1-\left(\dfrac{3}{5}\right)^{36}\bigg)\right)<0.
\end{equation}
Inequality \eqref{q003} is equivalent to the inequality
\begin{equation*}
\frac{\left(\dfrac{9}{25}\right)^{n-\sqrt{n}}}{1-\left(\dfrac{3}{5}\right)^{36}}<\frac{7}{37n^2},
\end{equation*}
and hence, for $q\in(0,\dfrac{9}{25}]$, we have
\begin{equation*}
\frac{q^{n-\sqrt{n}}}{1-q^{2n}}<\frac{\left(\dfrac{9}{25}\right)^{n-\sqrt{n}}}{1-\left(\dfrac{3}{5}\right)^{36}}<\frac{7}{37n^2}.
\end{equation*}

Thus, to complete the proof of the theorem, it remains to prove that for $n\geq9$ and ${q\in(\dfrac{9}{25},1)}$ the following implication is valid:
\begin{equation}\label{imp1}
\eqref{umova_n_0}\Rightarrow\eqref{umova_z}.
\end{equation}

Since
\begin{equation*}
\dfrac{1}{2}+\dfrac{2q}{(1+q^2)(1-q)}\leq \dfrac{1+q}{1-q},
\end{equation*}
relation \eqref{umova_n_0} implies the inequality
 \begin{equation*}
\dfrac{160}{57(n-\sqrt{n})}\; \dfrac{q}{(1-q)^2}<\left(\dfrac{1-q}{1+q}\right)^{\frac {4}{1-q^2}-1},
\end{equation*}
and hence, it also implies the equivalent inequality
\begin{equation}\label{sqrt_n0}
n-\sqrt{n}-\dfrac{160 q}{57(1-q)^2} \left(\dfrac{1+q}{1-q}\right)^{\frac {4}{1-q^2}-1}>0.
\end{equation}

Relation \eqref{sqrt_n0} implies
\begin{equation}\label{n1}
n>\dfrac{160q}{57(1-q)^2} \left(\dfrac{1+q}{1-q}\right)^{3}.
\end{equation}
Hence, for $n\geq9$ and ${q\in(0,1)}$ we have
\begin{equation}\label{imp2}
\eqref{umova_n_0}\Rightarrow\eqref{n1}.
\end{equation}

Further, let us prove that inequality \eqref{umova_z} for $n\geq9$ and ${q\in(0,1)}$ is implied by the inequality
\begin{equation}\label{n2}
n>\left(\frac{9(1+q)}{4(1-q)}\right)^{2}.
\end{equation}

Since for any $q\in(0,1)$ we have
\begin{equation*}
\ln \frac{1}{q}=2\sum\limits_{k=1}^\infty \frac{1}{2k-1}\left(\frac{1-q}{1+q}\right)^{2k-1}>2\frac{1-q}{1+q},
\end{equation*}
(see e.g. \cite[p.~53]{Gradshteyn_2007}), we have also
\begin{equation}\label{eq:104}
\left(\frac{9(1+q)}{4(1-q)}\right)^{2}>\left(\frac{9}{4\frac{1-q}{1+q}}\right)^{\frac {125}{79}}>\left(\frac{9}{2\ln 1/q}\right)^{\frac {125}{79}}.
\end{equation}
Relations \eqref{n2} and \eqref{eq:104} imply the inequality
\begin{equation*}
n>\left(\frac{9}{2\ln 1/q}\right)^{\frac {125}{79}},
\end{equation*}
which is equivalent to the following one:
\begin{equation}\label{eq:99}
\frac 23 n\ln \frac 1q>3 n^{\frac{46}{125}}.
\end{equation}
Since $\ln n<n^{\frac {46}{125}}$ for $n\in\mathbb{N}$ and $1-\dfrac{1}{\sqrt{n}}\geq\dfrac{2}{3}$ for $n\geq9$, relation~\eqref{eq:99} implies
\begin{equation}\label{eq:101}
n\left(1-\frac{1}{\sqrt{n}}\right)\ln \frac 1q>3\ln n.
\end{equation}
For $n\geq9$, from \eqref{eq:101} we obtain
\begin{equation*}
\dfrac{1}{q^{n}}> \dfrac{n^3}{q^{\sqrt{n}}}>\dfrac{9n^2}{q^{\sqrt{n}}}>\dfrac{38n^2}{7q^{\sqrt{n}}}=\dfrac{37n^2}{7q^{\sqrt{n}}}+\dfrac{n^2}{7q^{\sqrt{n}}}>\dfrac{37n^2}{7q^{\sqrt{n}}}+q^{n}.
\end{equation*}
Hence, for $n\geq9$ and $q\in(0, 1)$, we have
\begin{equation}\label{imp3}
\eqref{n2}\Rightarrow\eqref{umova_z}.
\end{equation}

It remains to prove that for ${q\in(\dfrac{9}{25},1)}$ and $n\geq9$ we have
\begin{equation}\label{imp4}
\eqref{n1}\Rightarrow\eqref{n2}.
\end{equation}
To this end, consider the difference of right-hand sides of inequalities \eqref{n1} and \eqref{n2} by setting
\begin{equation*}
v(q)=\dfrac{160q}{57(1-q)^2} \left(\dfrac{1+q}{1-q}\right)^{3}-\left(\frac{9(1+q)}{4(1-q)}\right)^{2}=
\end{equation*}
\begin{equation}\label{yq}
=
\left(\dfrac{1+q}{1-q}\right)^{2}\left(\dfrac{160q(1+q)}{57(1-q)^3} -\left(\frac{9}{4}\right)^{2}\right).
\end{equation}
Since ${q\in(\dfrac{9}{25},1)}$, we obtain
\begin{equation}\label{yq1}
\dfrac{160q(1+q)}{57(1-q)^3} -\left(\frac{9}{4}\right)^{2} >0.
\end{equation}
Relations \eqref{yq} and \eqref{yq1} imply the inequality $v(q)>0$ and also implication~\eqref{imp4}.
For ${q\in(\dfrac{9}{25},1)}$, relations \eqref{imp2}, \eqref{imp3}, and \eqref{imp4} imply \eqref{imp1}.
The theorem is proved.

As a consequence of proven theorem, we can formulate the following one.

\textbf{Theorem 3.}
\textit{Let $q\in(0,1)$.
Then for any $\beta\in \mathbb{R}$ and all numbers $n\geq n_q$ (where $n_q$ is the smallest number $n\geq9$, for which condition \eqref{umova_n_0} is satisfied) the following equalities take place}:
\begin{equation*}
d_{2n}(C_{\beta,\infty}^{q},C)=d_{2n-1}(C_{\beta,\infty}^q,C)= d_{2n-1}(C_{\beta,1}^q,L)=E_n(C_{\beta,\infty}^{q})_C=E_n(C_{\beta,1}^{q})_L=
\end{equation*}
\begin{equation}\label{bR}
=\|P_{q,\beta}\ast\varphi_n\|_C=\dfrac{4}{\pi}\left|\sum\limits_{\nu=0}^{\infty}\dfrac{q^{(2\nu+1)n}}{2\nu+1}\sin\left((2\nu+1)\theta_n\pi-\dfrac{\beta\pi}{2}\right)\right|,
\end{equation}
\textit{where $\theta_n=\theta_n(q,\beta)$ is the unique root of equation \eqref{theta} on $[0,1)$}.

\textit{In particular, whenever $n\geq n_q$ and $\beta\in \mathbb{Z}$, the following inequalities are satisfied}:
\begin{equation*}
d_{2n}(C_{\beta,\infty}^{q},C)=d_{2n-1}(C_{\beta,\infty}^q,C)= d_{2n-1}(C_{\beta,1}^q,L)=E_n(C_{\beta,\infty}^{q})_C=
\end{equation*}
\begin{equation}\label{b=0}
=E_n(C_{\beta,1}^{q})_L=\dfrac{4}{\pi}\arctan q^n, \; \beta=2k, \; k\in\mathbb{Z};
\end{equation}
\begin{equation*}
d_{2n}(C_{\beta,\infty}^{q},C)=d_{2n-1}(C_{\beta,\infty}^q,C)= d_{2n-1}(C_{\beta,1}^q,L)=E_n(C_{\beta,\infty}^{q})_C=
\end{equation*}
\begin{equation}\label{b=1}
=E_n(C_{\beta,1}^{q})_L=\dfrac{2}{\pi}\ln\dfrac{1+q^n}{1-q^n}, \; \beta=2k-1, \; k\in\mathbb{Z}.
\end{equation}

\textbf{Proof.}
In accordance with Theorem 2, for $n\geq n_q$ we have inequalities \eqref{dno1} and \eqref{dno2}. Combining those inequalities with formulas \eqref{E_n_rivnosti}, \eqref{d_m_E_n} and the relation $d_{2n-1}(C_{\beta,\infty}^{q},C)\geq d_{2n}(C_{\beta,\infty}^{q},C)$, we prove \eqref{bR} for $\beta\in \mathbb{R}$.
Now let us prove \eqref{b=0} and \eqref{b=1}.

For $\beta=2k$, $k\in\mathbb{Z}$, relation \eqref{theta_n} implies that the unique root of equation \eqref{theta} on $[0,1)$ is the value $\theta_n=\dfrac 12$.
Since in this case we have
\begin{equation*}
\sum\limits_{\nu=0}^{\infty}\dfrac{q^{(2\nu+1)n}}{2\nu+1}\sin\left((2\nu+1)\theta_n\pi-\dfrac{\beta\pi}{2}\right)=
\end{equation*}
\begin{equation*}
=(-1)^k\sum\limits_{\nu=0}^{\infty}(-1)^\nu\dfrac{q^{(2\nu+1)n}}{2\nu+1}=(-1)^k\arctan q^n,
\end{equation*}
relations \eqref{bR} imply \eqref{b=0}.

For $\beta=2k-1$, $k\in\mathbb{Z}$, relation \eqref{theta_n} implies that the root of equation \eqref{theta} is the value $\theta_n=0$.
Since
\begin{equation*}
\sum\limits_{\nu=0}^{\infty}\dfrac{q^{(2\nu+1)n}}{2\nu+1}\sin\left((2\nu+1)\theta_n\pi-\dfrac{\beta\pi}{2}\right)=
\end{equation*}
\begin{equation*}
=(-1)^{k}\sum\limits_{\nu=0}^{\infty}\dfrac{q^{(2\nu+1)n}}{2\nu+1}=(-1)^{k}\frac{1}{2}\ln\dfrac{1+q^n}{1-q^n},
\end{equation*}
relations \eqref{bR} imply equalities \eqref{b=1}.
The theorem is proved.

Note that Theorem 3 complements the results of \cite{Shevaldin_1992} in the case where $\dfrac{1}{5}<q<1$ and extends the mentioned result of \cite{Nguen_1994} onto the case where $\beta\in\mathbb{R}$.


\bibliographystyle{model1b-num-names}

\begin{thebibliography}{00}
 
 \bibitem {Bushanskiy_1978}
 A.V. Bushanskii, On the best harmonic approximation in the mean of some functions, in: Investigations in the Theory of Approximation of Functions and Their Applications, Institute of Mathematics, Ukrainian Academy of Sciences, Kiev, 1978, pp.~29--37 (in Russian).

\bibitem {Gradshteyn_2007}
I.S. Gradshtein, I.M. Ryzhik, Table of Integrals, Series, and Products, Academic Press, 2007.

\bibitem {Kornejchuk_1987}
 N.P. Korneichuk, Exact Constants in Approximation Theory, Nauka, Moscow, 1987 (in Russian).

\bibitem {Krein_1938}
M.G. Krein, On the theory of the best approximation of periodic functions, Dokl. Akad. Nauk SSSR 18 (4-5) (1938) 245--249 (in Russian).

\bibitem {Kuspel_1989}
A.K. Kushpel', Estimates of the diameters of convolution classes in the spaces $C$ and $L$, Ukr. Mat. Zh. 41 (8) (1989) 1070--1076 (in Russian).

\bibitem {Kuspel_1988}
A.K. Kushpel', Sharp estimates for the widths of convolution classes, Izv. Akad. Nauk SSSR, Ser. Mat. 52 (6) (1988) 1305--1322 (in Russian).

\bibitem {Nguen_1994}
Nguyen Thi Thieu Hoa, Extremal Problems on Some Classes of Smooth Periodic Functions, Doctoral-Degree Thesis (Physics and Mathematics), MIAN, Moscow, 1994 (in Russian).

\bibitem {Nikolskiy_1946}
S.M. Nikol'skii, Approximation of functions by trigonometric polynomials in the mean, Izv. Akad. Nauk SSSR Ser. Mat. 10 (1946) 207--256 (in Russian).

\bibitem {Pinkus_1985}
A. Pinkus, $n$-Widths in Approximation Theory, Springer-Verlag, Berlin, 1985.

\bibitem {Serdyuk_1998}
A.S. Serdyuk, Estimates for widths and best approximations of classes of convolutions of periodic functions, in: Fourier Series: Theory and Applications, Institute of Mathematics, Ukrainian Academy of Sciences, Kyiv, 1998, pp. 286--299 (in Ukrainian).

\bibitem {Serdyuk_1999_4}
A.S. Serdyuk, On the existence and uniqueness of a solution of the problem of uniform $SK$-spline-interpolation, Ukr. Mat. Zh. 51 (4) (1999) 486--492 (in Ukrainian).

\bibitem {Serdyuk_1999}
A.S. Serdyuk, Widths and best approximations for classes of convolutions of periodic functions, Ukr. Mat. Zh. 51 (5) (1999) 674--687 (in Ukrainian).

\bibitem {Serdyuk_2011}
 A.S. Serdyuk, I.V. Sokolenko, Asymptotic behavior of best approximations of classes of Poisson integrals of functions from $H_\omega$, Journal of Approximation Theory 163 (11) (2011) 1692--1706.

\bibitem {Shevaldin_1993}
V.T. Shevaldin, Lower bounds for the widths of classes of periodic functions with a bounded fractional derivative, Mat. Zametki 53 (2) (1993) 145--151 (in Russian).

\bibitem {Shevaldin_1989}
V.T. Shevaldin, Lower bounds on widths of classes of sourcewise representable functions, Trudy Mat. Inst. Steklov. 189 (1989) 185--200 (in Russian).

\bibitem {Shevaldin_1994}
V.T. Shevaldin, Lower estimates of the widths of classes of functions defined by a modulus of continuity, Izv. RAN. Ser. Mat. 58 (5) (1994) 172--188 (in Russian).

\bibitem {Shevaldin_1992}
V.T. Shevaldin, Widths of classes of convolutions with Poisson kernel, Mat. Zametki 51 (6) (1992) 126--136 (in Russian).

\bibitem {Stepanets_2005}
A.I. Stepanets, Methods of Approximation Theory, VSP, Leiden, 2005.

\bibitem {Serdyuk_1995}
 A.I. Stepanets, A.S. Serdyuk, Lower bounds for widths of classes of convolutions of periodic functions in the metrics of $C$ and $L$, Ukr. Mat. Zh. 47 (8) (1995) 1112--1121 (in Russian).

\bibitem {Stepanets_1994}
 A.I. Stepanets, A.S. Serdyuk, On the existence of interpolation $SK$-splines, Ukrain. Mat. Zh. 46 (11) (1994) 1546--1554 (in Russian).

\bibitem {Tihomirov_1976}
 V.M. Tikhomirov, Some Problems in Approximation Theory, Izd.~MGU, Moscow, 1976 (in Russian).



\end{thebibliography}



\end{document}